 \documentclass{ajour}

\newtheorem{lemmab}[theorem]{Lemma}
\newtheorem{propb}[theorem]{Proposition}
\newtheorem{corb}[theorem]{Corollary}

\newtheorem{defb}[theorem]{Definition}

\newtheorem{openquestion}[theorem]{Question}

\begin{document}

%


\authorrunninghead{Arturo Magidin}
\titlerunninghead{Amalgamation bases for nil-2 groups of odd exponent}





\title{Amalgamation bases for nil-2 groups of odd exponent%
\thanks*{2000 MSC: Primary 20E06,
20F18, 20D15, 08B25; Secondary 20E10.}%
\thanks*{The author was supported in part by CONACyT grant I~29922-E.}%
}


\author{Arturo Magidin}
\affil{Instituto de Matem\'aticas,
       Universidad Nacional Aut\'onoma de M\'exico\\
       Circuito Exterior, Ciudad Universitaria
       04510 Mexico City, Mexico}

\email{magidin@matem.unam.mx}

\abstract{We study amalgams and the strong, weak, and special
amalgamation bases in the varieties of nilpotent groups of class two
and exponent $n$, where $n$ is odd.  The main result is the
characterization of the special amalgamation bases for these
varieties. We also characterize the weak and strong bases.
For special amalgamation bases we show that there are groups
which are special bases in varieties of finite exponent but not in the
variety of all nil-2 groups, whereas for weak and strong bases we show
this is not the case. We also show that in these varieties, as well as
the variety of all nil-2 groups, a group has an absolute closure (in
the sense of Isbell) if and only if it is already absolutely closed,
i.e. if and only if it is a special amalgamation base.  }


\begin{article}


\zerosection{Introduction}

The main result in this paper is the characterization of the special
amalgamation bases in the variety ${\cal N}_2\cap{\cal B}_n$ of
nilpotent groups of class at most two and exponent $n$, with $n$
odd. We do this by first reducing to the case where $n$ is a prime
power, and then solving that case.  We then also extend the results to
characterize the strong and weak amalgamation bases. We also use the
characterization of the special amalgamation bases to show that in
these varieties, as well as the variety ${\cal N}_2$ of nilpotent
groups of class at most two, the only groups with an absolute closure
are the absolutely closed groups (we will recall the definitions
below).

In this section we will recall the basic definitions, and in Section~1
we will prove some technical lemmas. Since the case of special
amalgams turns out to be significantly more complicated than the
general case, we will start by studying them. In Section~2 we will
look at special amalgamation bases. In Section~3 we will use this
characterization to establish that there are special amalgamation
bases in the variety ${\cal N}_2\cap{\cal B}_{p^n}$ which are not
special amalgamation bases in ${\cal N}_2\cap{\cal B}_{p^{n+1}}$. We
will also look at absolute closures in these varieties, and address a
question of Higgins about absolutely closed algebras. In Section~4 we
will look at the strong and weak amalgamation bases. The
characterization will yield that in fact the only strong and weak
amalgamation bases in ${\cal N}_2\cap {\cal B}_{p^n}$ are those which
are amalgamation bases in ${\cal N}_2$. Finally, in Section 5, we will
discuss some related questions and some of the difficulties in the
case of $n$ even, as well as some partial results in that case.

Groups will be written multiplicatively, unless otherwise
specified. We will use $Z$ to denote the infinite cyclic group, which
we also write multiplicatively.  All maps are assumed to be group
morphisms unless we explicitly note otherwise. The multiplicative
identity of a group $G$ will be denoted by $e$, and we will use $e_G$
if there is danger of ambiguity. For a group $G$ and elements $x,y\in
G$, the commutator of $x$ and $y$ is $[x,y]=x^{-1}y^{-1}xy$; note that
$[x,y]^{-1}=[y,x]$. Given subsets $A,B$ of $G$, not necessarily
subgroups, $[A,B]$ denotes the subgroup of $G$ generated by all
commutators $[a,b]$ with $a\in A$ and $b\in B$. The commutator
subgroup of $G$ is the subgroup $[G,G]$, which is sometimes
denoted by $G'$. The center of $G$ is denoted by $Z(G)$. All
commutators will be written left-normed, so $[x,y,z]=[[x,y],z]$.

Let ${\cal N}_2$ denote the variety of all nilpotent groups of class
at most two; that is, groups $G$ such that $[G,G]\subseteq Z(G)$, or
equivalently, for which the identity $[x,y,z]=e$ holds. It is easy
to verify that for any nil-2 group (i{.}e{.} any nilpotent group of
class at most two) the following identities hold, so we will use them
without comment throughout the paper.

\begin{propb}
Let $G\in {\cal N}_2$. For all $x,y,z\in G$ and all integers $n$:
\begin{itemize}
\item[(a)]\  $[xy,z]=[x,z][y,z]$; $[x,yz]=[x,y][x,z]$.
\item[(b)]\  $[x^n,y]=[x,y]^n=[x,y^n]$.
\item[(c)]\  $(xy)^n = x^n y^n[y,x]^{n(n-1)/2}$.
\item[(d)]\ The value of $[x,y]$ depends only on the congruence
classes of $x$ and~$y$ modulo $G'$ (in fact, modulo $Z(G)$). 
\end{itemize}
\end{propb}

From these properties, it is easy to see that if $G\in {\cal N}_2$ is
generated by elements of exponent $n$ with $n$ odd, then $G$ itself is of exponent
$n$, and, if $n$ is even, of exponent $2n$. Also, if $K$
is generated by a subgroup $G$ of exponent $2n$ and elements of
exponent $n$, then $K$ is also of exponent $2n$, whether $n$ is even
or not.

We use ${\cal B}_n$ to denote the variety of all groups of exponent
$n>0$. Thus, the variety of all nilpotent groups of class at most two
and exponent $n$ is denoted by~${\cal N}_2\cap{\cal B}_n$.

For the remainder of this paper, all groups will be assumed to lie in
${\cal N}_2$ unless otherwise specified. Any presentation of a group
will be understood to be a presentation in ${\cal N}_2$; that is, the
identities on ${\cal N}_2$ will be imposed on the group, as well as
all the relations specified in the presentation. 

Given $A,B\in{\cal N}_2$, every element of their coproduct
$A\coprod^{{\cal N}_2} B$ has a unique expression of the form
$\alpha\beta\gamma$, where $\alpha\in A$, $\beta\in B$, and $\gamma\in
[A,B]$. A theorem of T.~MacHenry \cite{machenry} states that the subgroup
$[A,B]$ of $A\coprod^{{\cal N}_2}B$ is isomorphic to the tensor
product $A^{\rm ab}\otimes B^{\rm ab}$. 

Recall that an ${\cal N}_2$-amalgam of two groups $A,C\in {\cal N}_2$
with core $B$ consists of groups $A$, $B$, and $C$, equipped with one
to one group morphisms
\begin{eqnarray*}
\Phi_A\colon B & \to & A\\
\Phi_C\colon B & \to & C.
\end{eqnarray*}

To simplify notation, we denote this situation by $(A,C;B)$. To say
the amalgam $(A,C;B)$ is \textit{(weakly) embeddable in ${\cal N}_2$}
means that there exists a group $M$ in ${\cal N}_2$ and one to one
group morphisms
\[\lambda_A\colon A\to M,\qquad \lambda_C\colon C\to M, \qquad
\lambda\colon B\to M,\]
such that
\[\lambda_A\circ\Phi_A = \lambda\qquad{\rm and}\qquad \lambda_C\circ\Phi_C=\lambda.\]

When we examine whether or not the amalgam $(A,C;B)$ is embeddable,
the obvious candidate for $M$ is the coproduct with amalgamation of
$A$ and $C$ over $B$, denoted by $A\coprod^{{\cal N}_2}_{B} C$. This
coproduct is sometimes called the ${\cal N}_2$-free product with
amalgamation. We say that $(A,C;B)$ is \textit{weakly embeddable} (in
${\cal N}_2$) if no two distinct elements of $A$ are identified with
each other in the coproduct with amalgamation, and similarly with two
distinct elements of $C$. Note that weak embeddability does not
preclude the possibility that an element $x$ of $A\setminus B$ be
identified with an element $y$ of $C\setminus B$ in
$A\coprod_{B}^{{\cal N}_2} C$. We say that $(A,B;C)$ is
\textit{strongly embeddable} (in ${\cal N}_2$) if there is also no
identification between elements of $A\setminus B$ and elements of
$C\setminus B$. By \textit{special amalgam} we mean an amalgam
$(A,A';B)$, where there is an isomorphism $\psi\colon A\to A'$ such
that $\psi\circ\Phi_A=\Phi_{A'}$. In this case, we usually write
$(A,A;B)$ with $\psi={\rm id}_A$ being understood. Since special
amalgams are always weakly embeddable (say, by mapping to $A$), with
special amalgams we are mostly interested in whether or not they can
be strongly embedded.

Finally, recall that a group $G$, which belongs to a class of groups
${\cal C}$, is said to be a \textit{weak amalgamation base} for ${\cal
C}$ if and only if every amalgam of ${\cal C}$-groups with core $G$ is
weakly embeddable in ${\cal C}$; it is a \textit{strong amalgamation
base} if every such amalgam is strongly embeddable; and it is a
\textit{special amalgamation base} if every special amalgam of ${\cal
C}$-groups with core $G$ is strongly embeddable.

Special amalgams are closely related to the concept of
\textit{dominion}. Recall that Isbell \cite{isbellone} defines for a
category $\cal C$ of algebras (in the sense of Universal Algebra) of a
fixed type $\Omega$, and an algebra $A\in {\cal C}$ and subalgebra $B$
of $A$, the \textit{dominion of $B$ in $A$} (in the category $\cal C$)
to be the intersection of all equalizers containing $B$. Explicitly,
\[{\rm dom}_A^{{\cal C}}(B) = \left\{a\in A | \forall f,g\colon A\to C,
\ {\rm if}\ f|_B=g|_B\ {\rm then}\ f(a)=g(a)\right\}\]
where $C$ ranges over all algebras $C\in {\cal C}$, and $f,g$ are
morphisms. The connection between amalgams and dominions when
working in a variety is through special amalgams:
letting $A'$ be an isomorphic copy of
$A$, and $M=A\coprod^{{\cal C}}_B A'$, we have that
\[{\rm dom}_A^{{\cal C}}(B) = A\cap A'\subseteq M\]
where we have identified $B$ with its common image in $A$ and $A'$. In
other words, ${\rm dom}_A^{{\cal C}}(B)$ is the smallest subalgebra
$D$ of~$A$ such that $B\subseteq D$ and the amalgam $(A,A;D)$ is
strongly embeddable. If ${\rm dom}_A^{{\cal C}}(B)=B$, we say that the
dominion of~$B$ is ``trivial'' (meaning it is as small as possible),
and we say it is ``nontrivial'' otherwise.

In general, ${\rm dom}_A^{{\cal C}}(-)$ is a closure operator on the
lattice of subalgebras of~$A$. If we are working in a variety of
groups (i{.}e{.}, a full subcategory of ${\cal G}roup$ which is closed
under taking subgroups, quotients, and arbitrary direct products),
then the dominion construction respects finite direct products (that
is, if $H_1<G_1$ and $H_2<G_2$, then the dominion of $H_1\times H_2$
in $G_1\times G_2$ is the product of the dominions of $H_1$ in $G_1$
and of $H_2$ in $G_2$), and also respects quotients: if
$N\triangleleft G$ is contained in $H$, then
\[{\rm dom}_{G/N}^{\cal V}(H/N) = \left({\rm dom}_G^{\cal
V}(H)\right)\bigm/ N.\]
For a proof of these assertions we direct the reader to
\cite{nildoms}.

Because of the connection with dominions, special amalgamation bases
are also said to be \textit{absolutely closed.} We direct the reader
to the survey paper by Higgins~\cite{episandamalgs} for a more
complete discussion of amalgams and their connection with dominions.

\section{Preliminary Results}

In the variety ${\cal N}_2$ not every amalgam is weakly
embeddable, and not every dominion is trivial; see for example
\cite{nildoms} and \cite{nonembed}. 

The following easy observations will be key:

\begin{propb}
\label{bothfiniteexp}
Let $(G,K;H)$ be an amalgam of ${\cal N}_2$-groups. Assume both $G$ and
$K$ are of exponent $n$. If $n$ is odd, then $(G,K;H)$ is weakly
(resp.\ strongly) embeddable into an ${\cal N}_2$-group if and only if
$(G,K;H)$ is weakly (resp.\ strongly) embeddable into an ${\cal
N}_2\cap{\cal B}_n$-group. If $n$ is even, then $(G,K;H)$ is weakly
(resp.\ strongly) embeddable into an ${\cal N}_2$-group if and only if
$(G,K;H)$ is weakly (resp.\ strongly) embeddable into an ${\cal
N}_2\cap{\cal B}_{2n}$-group.
\end{propb}

\begin{proof}
Clearly, we may assume that a group into which we embedd the amalgam
$(G,K;H)$ is generated by $G$ and $K$. For $n$ odd this means the
group itself is of exponent $n$, and for $n$ even of exponent $2n$.
\end{proof}

\begin{propb}
\label{pparts}
Let $(G,K;H)$ be an amalgam of ${\cal N}_2$-groups. If both $G$ and
$K$ are of exponent $n$, then $(G,K;H)$ is weakly (resp.\ strongly)
embeddable into an ${\cal N}_2$-group if and only if for every prime
$p$ dividing $n$, the amalgam $(G_p,K_p;H_p)$ is weakly (resp.\
strongly) embeddable, where $G_p$ is the $p$-part of $G$, and likewise
for $K_p$ and $H_p$.
\end{propb}

\begin{proof}
Since $G$, $K$, and $H$ are of finite exponent, they are the direct
product of their $p$-parts, i.e. $G=\prod G_p$, $K=\prod K_p$, and
$H=\prod H_p$. Clearly, an embedding of $(G,K;H)$ provides embeddings
for $(G_p,K_p;H_p)$ for each $p$; conversely, if we have embeddings
into $M_p$ for each $(G_p,K_p;H_p)$, then $\prod M_p$ will give an
embedding for $(G,K;H)$.
\end{proof}

In view of the above two propositions, we may restrict ourselves to
the varieties ${\cal N}_2\cap {\cal B}_{p^n}$ with $p$ a prime.  In
this paper we will further restrict ourselves to the case of $p$ an
odd prime, except for some comments near the end.

Note that the varieties ${\cal N}_2\cap{\cal B}_{p^n}$ are smaller
than the variety ${\cal N}_2$. Clearly if a group $G\in {\cal
N}_2\cap{\cal B}_{p^n}$ is an amalgamation base for ${\cal N}_2$
(weak, strong, or special), then it is also an amalgamation base (of
the same type) for ${\cal N}_2\cap{\cal B}_{p^n}$, by
Proposition~\ref{bothfiniteexp}; the converse, however, is not immediate and
could fail. So a characterization of such bases is an interesting
problem, whether it turns out that the converse does hold or that it
does not.

It will be helpful to recall a theorem on adjunction of roots to
${\cal N}_2$.

\begin{theorem}[Saracino, Theorem 2.1 in \cite{saracino}]
\label{genrootadjunction}
Let $G$ be a nilpotent group of class at most two, let $m>0$, let
$\mathbf{n}$ be an $m$-tuple of positive integers, and let $\mathbf{g}$
be an $m$-tuple of elements of~$G$. Then there exists a nilpotent
group $K$ of class two, containing $G$, and which contains an $n_i$-th
root for $g_i$ ($1\leq i\leq m$) if and only if for every $m\times m$
array $\{c_{ij}\}$ of integers such that $n_i c_{ij}=n_j c_{ji}$ for
all $i$ and $j$, and for all $y_1,\ldots, y_m\in G$,
\[ {\it if}\ y_j^{n_j}\equiv
\prod_{i=1}^{m}g_i^{c_{ij}}\pmod{G'}\qquad {\it then}\
\prod_{j=1}^{m} [y_j,g_j] = e.\]
\end{theorem}

Note that the theorem implies that we can always adjoin $n_i$-th roots
to a family of central elements, and so to a family of
commutators. Therefore, if the $g_i$ are $n_i$-th powers modulo the
commutator in $G$, then there is an extension with $n_i$-th roots for
the $g_i$.

We also have the following two results:

\begin{lemmab}
\label{centralroots}
Let $G \in {\cal N}_2$ and let $g_1,\ldots,g_r\in Z(G)$. For all
$r$-tuple of postive integers $\mathbf{n}=(n_1,\ldots,n_r)$ there is
an extension $K$ of $G$ with an $n_i$-th root $h_i$ of $g_i$, 
and such that $h_i$ is central in~$K$, for $1\leq i\leq r$.
\end{lemmab}
\begin{proof}
Let $a_i$ be the order of $g_i$ ($a_i=0$ if $g_i$ is not torsion);
then consider the group
$F = G \times (Z/n_1a_1Z)\times\cdots\times (Z/n_ra_rZ)$.
Denote the generator of the cyclic groups by $h_i$, and mod out by the
subgroup generated by $g_ih_i^{-n_i}$. Since the $g_i$ are central, the
subgroup is normal and this works out.
\end{proof}

\begin{lemmab}
\label{makecentralacomm}
Let $G\in {\cal N}_2$ and let $g_1,\ldots,g_r\in Z(G)$. Then there
exists an extension $K$ of $G$, with $K\in{\cal N}_2$, and such that
$g_i\in [K,K]$ for every $i$. If $G$ is of odd exponent $n$, then we
can choose $K$ of exponent $n$. If $G$ is of even exponent $2n$, then
we can always choose $K$ of exponent $4n$. If $G$ is of even exponent
$2n$ and the $g_i$ are of exponent $n$, we may in fact choose $K$ of
exponent $2n$.
\end{lemmab}

\begin{proof}
Let $a_i$ be the order of $g_i$, with $a_i=0$ if $g_i$ is not
torsion. Consider the group
\[F = G \times \left((Z/a_1Z)\coprod\nolimits^{{\cal
N}_2}(Z/a_1Z)\right)\times\cdots \times \left(
(Z/a_rZ)\coprod\nolimits^{{\cal N}_2}(Z/a_rZ)\right)\] and denote the
generators of the cyclic groups by $t_{i1}$ and $t_{i2}$. Then mod out
by the subgroup generated by $g_i[t_{i2},t_{i1}]$. Again, since the
$g_i$ are central, this subgroup is normal. If $G$ is of odd exponent
$n$, then $F$ is generated by elements of exponent $n$, and so it is
again of exponent $n$, hence so is $F/N$.

If $G$ is of exponent $2n$, then $F$ is generated by elements of
exponent $2n$ and so $F/N$ is of exponent at most $4n$. If, furthermore,
the $g_i$ are of exponent $n$, then $F$ is generated by a group of
exponent $2n$ and elements of exponent $n$, so $F$ is of exponent
$2n$, and then so is $F/N$.
\end{proof}

The following lemma will be used several times:

\begin{lemmab}
\label{keylemma}
Let $G\in {\cal N}_2$, $x,y\in G$. Let $K$ be any ${\cal
N}_2$-extension of~$G$, and assume that for some $n>0$, $r,s\in K$,
$r',s'\in K'$, we have that $x=r^nr'$ and $y=s^ns'$. For any integers
$a$, $b$, $c$, and $j$, and all $g_1,g_2\in G$,
\[ \hbox{{\it if}}\qquad
\begin{array}{rcl}
g_1^n & \equiv & x^a y^b \pmod{K'} \\
g_2^n & \equiv & x^{b+j} y^c \pmod{K'}
\end{array}
\quad\hbox{\it then}\quad
[r,s]^{jn} = [g_1,x][g_2,y].\]
In particular, $[r,s]^{jn}$ lies in~$G$.
\end{lemmab}

\begin{proof}
Since $xr^{-n}$ and $ys^{-n}$ are commutators in $K$, they are
central. Therefore:
\begin{eqnarray*}
e & = & [g_1r^{-a}s^{-b-j},xr^{-n}][g_2r^{-b}s^{-c},ys^{-n}] \\
  & = & [g_1,x][g_2,y][r^{-a}s^{-b-j},x] [r^{-b}s^{-c},y] [g_1,r^{-n}]\\
  &   &\quad        [g_2,s^{-n}][s^{-b-j},r^{-n}][r^{-b},s^{-n}]\\
  & = & [g_1,x][g_2,y][r,x^{-a}y^{-b}][s,x^{-b-j}y^{-c}][g_1^{-n},r]
  [g_2^{-n},s][r,s]^{-jn} \\
  & = &
  [g_1,x][g_2,y][r,g_1^{-n}][s,g_2^{-n}][g_1^{-n},r][g_2^{-n},s][r,s]^{-jn}\\
  & = & [g_1,x][g_2,y][r,s]^{-jn}.
\end{eqnarray*}
Therefore, $[r,s]^{jn}=[g_1,x][g_2,y]$, as claimed.
\end{proof}

We note that the condition depends only on the classes of $x$ and $y$
modulo $G'$; moreover, the congruences are actually symmetric on $x$
and $y$. If you exchange $x$ and $y$, then simply exchange $a$ with
$c$, and replace, for any $j\geq 0$, $b$ with $-b-j$, $g_1$ with
$g_2^{-1}$, and $g_2$ with $g_1^{-1}$.

\begin{remark} In many of our applications, the congruences in
Lemma \ref{keylemma} will hold modulo $G'$ rather than $K'$; however,
since $G'\subseteq K'$, the Lemma will apply in that case as well.
\end{remark}

\begin{remark} If $x$ or $y$ lie in $G^nG'$, then we can always find
integers $a$, $b$, and $c$, and elements $g_1$, $g_2\in G$ satisfying
the congruences in Lemma~\ref{keylemma}, so in that case
$[r_1,r_2]^{jn}\in G$. For example, if $x\in G^nG'$, write $x=g^ng'$
and let $a=b=c=0$, $g_1=e$, and $g_2=g^j$.
\end{remark}

\section{Special amalgamation bases}

In this section we characterize the special amalgamation bases in the
variety ${\cal N}_2\cap{\cal B}_{p^n}$, with $p$ an odd prime.

Note that if $K \in {\cal N}_2\cap{\cal B}_{p^n}$, with $p$ an odd
prime, then, in view of Proposition~\ref{bothfiniteexp}, the dominion
of a subgroup $G$ of~$K$ in ${\cal N}_2\cap{\cal B}_{p^n}$ is the same
as the dominion in ${\cal N}_2$.  The description of dominions is:

\begin{theorem}[Theorem 3.31 in~\cite{nildoms}]
\label{nildomsdesc}
Let $K\in{\cal N}_2$, $G$ a subgroup of~$K$. Let $D$ be the subgroup
of~$K$ generated by $G$ and all elements $[r,s]^q$, where $q>0$ and
$r^q,s^q\in G[K,K]$. Then $D={\rm dom}_K^{{\cal N}_2}(G)$.
\end{theorem}

\begin{remark}
It is an easy exercise to verify that we can restrict $q$ to prime
powers.
\end{remark}

\begin{remark} Also note that if $r^q,s^q\in G[K,K]$, and one of them
has a $q$-th root in~$G$ modulo $[K,K]$, then $[r,s]^q\in G$. For if
$r^q$ has such a $q$-th root, then 
$r^q\equiv g^q \pmod{K'}$ for some $g\in G$, and
\[[r,s]^q = [r^q,s] = [g^q,s] = [g,s^q]\]
which lies in $G$ since $s^q$ is congruent to an element of $G$ modulo
$K'$. A similar calculation holds if $s^q$ has a $q$-th root.
\end{remark}

In particular, if $K$ is of exponent $p^n$ with $p$ an odd prime, then in
Theorem~\ref{nildomsdesc} we can restrict $q$ to prime powers $p^i$, with
$1\leq i < n$.

In order to characterize special amalgamation bases, we will start with a
group $G$, and we will be interested to know when, for any overgroup
$K$, an element $[r,s]^q$ which lies in the dominion must also lie in
$G$. The idea is to look at pairs of elements of~$G$; call these $x$
and $y$; we are thinking that $x=r^qr'$ and $y=s^qs'$, with $r,s\in K$
and $r',s'\in K'$. We want $[r,s]^q\in G$. Now, for any pair of
elements, it is possible that no such extension exists (for example,
if no ${\cal N}_2$-extension exists with $q$-th roots for both $x$ and
$y$); or else that the situation is possible, in which case we want a
condition like that found in Lemma~\ref{keylemma} to ensure that
$[r,s]^q\in G$. So we will have a disjunction of conditions, where
some of them are used to rule out the possiblity that $x=r^qr'$ and
$y=s^qs'$ in some overgroup, and the rest are used to guarantee that,
even if we {\it can} find such $K$, $r$, $s$, $r'$, and $s'$, we will
nevertheless have $[r,s]^q\in G$.

When we work in ${\cal N}_2$, the only case in which there is no $K$
with elements $r$, $s$, $r'$, and $s'$ with the given properties is
when no ${\cal N}_2$-extension with $q$-th roots for $x$ and $y$
exists at all. This gives:

\begin{theorem}[Theorem 2.9 in~\cite{absclosed}]
\label{classiffullvar}
Let $G\in{\cal N}_2$. Then $G$ is absolutely closed in ${\cal N}_2$ if
and only if for all $x$, $y\in G$, and for all $n>0$, one of the
following holds:
\begin{itemize}
\item[(i)] There exist integers $a$, $b$, and~$c$, and elements
$g_1,g_2\in G$ such that
\begin{eqnarray*}
g_1^n & \equiv & x^ay^b\pmod{G'}\\
g_2^n & \equiv & x^by^c\pmod{G'}
\end{eqnarray*}
and $[g_1,x][g_2,y]\not=e$; or
\item[(ii)] There exist integers $a$, $b$, and~$c$, and elements
$g_1,g_2\in G$ such that
\begin{eqnarray*}
g_1^n & \equiv & x^ay^b\pmod{G'}\\
g_2^n & \equiv & x^{b+1}y^c\pmod{G'}.
\end{eqnarray*}
\end{itemize}
\end{theorem}

Our main result follows. Its statement is more complicated than
Theorem~\ref{classiffullvar} because the conditions for the existence
of an overgroup~$K$ in ${\cal N}_2\cap{\cal B}_{p^n}$ with the property
explained above are more complicated.

\begin{theorem}
\label{classifabsclosed}
Let $G\in{\cal N}_2\cap{\cal B}_{p^n}$, where $p$ is an odd prime and
$n\geq 1$. Then $G$ is absolutely closed in ${\cal N}_2\cap{\cal
B}_{p^n}$ if and only if for every $x,y\in G$, and each
each $i$, $1\leq i\leq n-1$, one of the following holds:
\begin{itemize}
\item[(a)]\ There exist $g_1,g_2\in G$, and integers $a,b,c$ such that
\begin{eqnarray*}
g_1^{p^i} & \equiv & x^a y^b \pmod{G'} \\
g_2^{p^i} & \equiv & x^b y^c \pmod{G'}
\end{eqnarray*}
and $[g_1,x][g_2,y]\not=e$; or

\item[(b)]\ $[x,y]^{p^{\alpha}}\not=e$ where $\alpha=\max\{0,n-2i\}$; or

\item[(c)]\ $2i>n$ and there exist $g_1,g_2\in G$, integers $a$, $b$,
and $c$, and an integer $j\in\{p^{n-i},2p^{n-i},\ldots,p^i\}$ such that
\begin{eqnarray*}
g_1^{p^i} & \equiv & x^a y^b \pmod{G'}\\
g_2^{p^i} & \equiv & x^{b+j} y^c \pmod{G'}
\end{eqnarray*}
and $[g_1,x][g_2,y]\not=e$; or

\item[(d)]\ There exist $g_1,g_2\in G$ and integers $a,b,c$ such that
\begin{eqnarray*}
g_1^{p^i} & \equiv & x^a y^b \pmod{G'}\\
g_2^{p^i} & \equiv & x^{b+1} y^c \pmod{G'}.
\end{eqnarray*}
\end{itemize}
\end{theorem}

\begin{remark} Note that the conditions on $x$ and $y$ depend only on
their congruence class modulo $G'$, and that they all are
symmetric on $x$ and $y$. Also note that in (a), (b), and~(c), we can
restrict $a$, $b$, and $c$ to being integers between $0$ and $p-1$.
\end{remark}

\begin{remark} One should think of conditions (a), (b), and (c) as
``pre-emptive,'' and condition (d) as being ``reactive.''
Condition (a) is nothing more than the statement that there is no
${\cal N}_2$-extension of $G$ with $p^{i}$-th roots for both $x$ and
$y$.  The way in which condition (a) works to prevent the existence of
an extension in which $x$ and $y$ are $p^i$-th powers modulo the
commutator is clear; for condition (b), note that if $x$ and $y$ are
both $p^i$-th powers modulo the commutator, then their commutator is
really the $p^{2i}$-th power of a commutator in the overgroup, so
$[x,y]^{p^{\alpha}}$ should be trivial. For $i$ bigger than $n/2$,
condition (c) ensures that $x^{n-i}$ would not commute with the
possible $p^i$-th root of $y$, which would be a problem if our groups
are really groups of exponent $p^n$ and $x$ is a $p^i$-th power modulo
the commutator. The conjunction of the negation of (a), (b), and (c)
are thus necessary for $x$ and $y$ to be $p^i$-th powers modulo the
commutator in some extension of exponent $p^n$, and the proof will
show that they are also sufficient; condition (d) comes to the rescue
in the situation when $x$ and $y$ \textit{are} $p^i$-th powers modulo
the commutator in some extension of exponent $p^n$.
\end{remark}

\begin{proof}
Suppose that $G$ satisfies the conditions, and let $K$ be any
overgroup of $G$ with $K\in{\cal N}_2\cap{\cal B}_{p^n}$. Let
$r,s\in K$ and $r',s'\in K'$ be such that $r^{p^i}r',s^{p^i}s'\in G$
for some $i$, $1\leq i\leq n-1$. Let $x=r^{p^i}r'$, $y=s^{p^i}s'$. We
want to verify that $[r,s]^{p^i}\in G$.

Note that since $x$ and $y$ are $p^i$-th powers modulo the commutator
in $K$, then there is an extension of $K$ in ${\cal N}_2$ with
$p^i$-th roots for both $x$ and $y$; therefore condition (a) cannot
hold for $x$, $y$, and $i$. Likewise, condition (b) does not hold: if
$2i\leq n$, then $[x,y]^{p^{n-2i}} = [r,s]^{p^n}=e$, since $K$ is of
exponent $p^n$; and if $2i>n$, then $e=[r,s]^{p^{2i}}=[x,y]$.

We also have that $x$, $y$, and $i$ cannot satisfy (c). For suppose
that $2i>n$ and for some $j\in\{p^{n-i},2p^{n-i},\ldots,p^i\}$, we 
have the congruences
\begin{eqnarray*}
g_1^{p^i} & \equiv & x^a y^b\pmod{G'}\\
g_2^{p^i} & \equiv & x^{b+j} y^c\pmod{G'}
\end{eqnarray*}
with $[g_1,x][g_2,y]\not=e$. Then by Lemma~\ref{keylemma}, we have
$[r,s]^{jp^i}\not=e$. However, for all such $j$, $jp^i$ is a multiple
of $p^n$, hence $[r,s]^{jp^i}=e$ in $K$. Since we fail conditions
(a)--(c), condition (d) must
be satisfied, and then by Lemma~\ref{keylemma} it follows that
\[[r,s]^{p^i} = [g_1,x][g_2,y]\in G.\]
So if the conditions hold, then $G$ is absolutely closed, as claimed.

Necessity is more complicated. We assume that $G$, $x$, $y$, and $i$
do not satisfy the conditions (a)--(d). We will proceed as follows:
first we construct an extension where $x$ and $y$ have $p^i$-th roots
$r$ and $s$, respectively, and where $[r,s]^{p^i}\notin G$. Then, if
necessary, we make $[r,s]^{p^n}$ trivial. Then we verify that
$x^{p^{n-i}}$ and $y^{p^{n-i}}$ are central in this extension, so we
can adjoin \textit{central} $p^n$-th roots $t$ and $v$, respectively,
to each. Next, we make $t^{p^i}$ and $v^{p^i}$ into commutators. Then
$rt^{-1}$ and $sv^{-1}$ are of exponent $p^n$, and their $p^i$-th
powers are, modulo commutators, $x$ and $y$. Finally, we select a
suitable subgroup which is generated by elements of exponent $p^n$. If
we ensure that throughout the process $[r,s]^{p^i}\notin G$, we will
have an example where $G$ is not equal to its dominion. We now proceed
with this program:

STEP 1. To facilitate notation, we write $x_1$ and $x_2$ rather than
$x$ and $y$ (respectively) in this step. Let $K_0=\bigl(G\coprod^{{\cal
N}_2}(Z\coprod^{{\cal N}_2}Z)\bigr)/N$, where we denote the generators of
the two copies of $Z$ by $r_1$ and $r_2$, and $N$ is the normal
subgroup generated by $x_1r_1^{-p^i}$, $x_2r_2^{-p^i}$.

We need to verify that $N\cap G=\{e\}$ and that for every $g\in G$ we have
$g[r_1,r_2]^{-p^i}\notin G$. The argument is the same as that found in the
proof of Theorem 2.9 in \cite{absclosed}, so we only sketch it here:

A general element of $N$ can be written as
\[\prod_{j=1}^2\left(\prod_{k=1}^{s_j} (b_{jk}z_{jk})(x_j
r_j^{-p^i})^{\epsilon_{jk}}(b_{jk}z_{jk})^{-1}\right)\]
where $b_{jk}\in G$, $s_j$ is a positive integer, $\epsilon_{jk}=\pm
1$, and $z_{jk}=r_1^{a_{jk1}}r_2^{a_{jk2}}$. 
By rewriting the element and expanding brackets, we get that the element equals:
\[\prod_{j=1}^2\left(\prod_{k=1}^{s_j}
[b_{jk},x_j]^{\epsilon_{jk}}[b_{jk},r_j^{-p^i}]^{\epsilon_{jk}}
[z_{jk},x_j]^{\epsilon_{jk}} [z_{jk},r_j^{-p^i}]^{\epsilon_{jk}}\right)
(x_jr_j^{-p^i})^{t_j}\]
where $t_j=\sum_{k=1}^{s_j}\epsilon_{jk}$.

Say that this is equal to an element of $G\coprod^{{\cal
N}_2}(Z\coprod^{{\cal N}_2}Z)$ of the form $g[r_1,r_2]^{qp^i}$ for
some integer $q$ and some $g\in G$.  Writing this as
$\alpha\beta\gamma$ with $\alpha\in G$, $\beta\in Z\coprod^{{\cal
N}_2} Z$, $\gamma\in [G,Z\coprod^{{\cal N}_2}Z]$, we obtain equations:
\begin{eqnarray*}
g & = & [h_1,x_1][h_2,x_2] \\
{[r_2,r_2]}^{qp^i} & = & [r_1,r_2]^{p^i(c_{21}-c_{12})}\\
e & = &
[h_1^{-p^i}x_1^{c_{11}}x_2^{c_{21}},r_1][h_2^{-p^i}x_1^{c_{12}}x_2^{c_{22}},
r_2]
\end{eqnarray*}
where $h_j=\prod_k b_{jk}^{\epsilon_{jk}}$, and $c_{ij}=-\sum\epsilon_{ik}a_{ikj}$.
This means that
\begin{eqnarray*}
h_1^{p^i} & \equiv & x_1^{c_{11}}x_2^{c_{21}}\pmod{G'}\\
h_2^{p^i} & \equiv & x_1^{c_{12}}x_2^{c_{22}}\pmod{G'}.
\end{eqnarray*}

For $q=0$, we have $c_{12}=c_{21}$, so the failure of condition (a)
shows that $N\cap G=\{e\}$.  By setting $q=-1$, we get that
$c_{12}=c_{21}+1$, so we are in the situation of (d); the fact that
(d) fails in $G$ means that no such solution exists, so no element of
the form $g[r_1,r_2]^{-p^i}$ lies in $N$. Therefore, in $K_0$ we have
$[r_1,r_2]^{p^i}$ lies in the dominion of $G$ but not in~$G$. We
change the names of $r_1$ and $r_2$ to $r$ and $s$ respectively, and
we again write $x$ and $y$ instead of $x_1$ and $x_2$, respectively,
to reduce the number of indices.

So we have that $K_0$ is an extension of $G$, generated by $G$, $r$
and $s$; that $r^{p^i}=x$, $s^{p^i}=y$, and that $[r,s]^{p^i}\notin
G$. We also note that $[r,s]^{p^{2i+\alpha}}\in N$, with $\alpha$ as
defined in (b). Indeed, $[xr^{-p^i},ys^{-p^i}]^{p^{\alpha}}\in N$, and
\[[xr^{-p^i},ys^{-p^i}]^{p^{\alpha}} = [x,y]^{p^{\alpha}}
[x,s]^{-p^{i+\alpha}} [r,y]^{-p^{i+\alpha}} [r,s]^{p^{2i+\alpha}}.\]
However, $[x,y]^{p^{\alpha}}=e$ since $G$ fails (b). And also
\[[x,s]^{-p^{i+\alpha}} =[x,y]^{p^{\alpha}}[x,s^{-p^i}]^{p^{\alpha}}
= [x,ys^{-p^i}]^{p^{\alpha}} = (ys^{-p^i})^x ys^{-p^i}\in N,\]
and similarly with $[r,y]^{-p^{i+\alpha}}$. So
$[r,s]^{p^{2i+\alpha}}\in N$.

We also claim that if $2i>n$ and $g[r,s]^{kp^n}\in N$ with
$k\in\{1,2,\ldots,p^{2i-n}\}$, then $g=e$. Indeed, if we write
\[g[r,s]^{kp^n} = g[r,s]^{p^i(kp^{n-i})}\]
we note that $kp^{n-i}\in\{p^{n-i},2p^{n-i},\ldots,p^i\}$, and so by
construction of $K_0$ we must have $g_1,g_2\in G$, and $a,b,c$
integers with
\begin{eqnarray*}
g_1 & \equiv & x^a y^b\pmod{G'}\\
g_2 & \equiv & x^{b+kp^{n-i}} y^c \pmod{G'}
\end{eqnarray*}
with $g=[g_1,x][g_2,y]$. However, the failure of condition (c) means
that in this case, $g=e$, as claimed.

STEP 2. Let $K_1 = K_0/M$ where $M$ is the normal subgroup generated
by $[r,s]^{p^n}$.  We claim that $M\cap G= \{e\}$, and that
$g[r,s]^{p^i}\notin M$ for each $g\in G$. 

Indeed, if $2i\leq n$, then $M$ is actually trivial already, since
$[r,s]^{p^{2i+\alpha}}\in N$; if $2i>n$, then we know that $[r,s]^{kp^n}\notin
G$ for $k=1,2,\ldots,p^{2i-n}$ again by the observations above; thus
$M\cap G = \{e\}$; and we can verify by inspection that
$g[r,s]^{p^i}\notin M$ for each $g\in G$, since $[r,s]^{p^n}$ is
central. 

So $[r,s]^{p^i}$ is still in the dominion of $G$ but not in $G$.

STEP 3. We note that $x^{p^{n-i}}$ and $y^{p^{n-i}}$ are central in
$K_1$. That they are central in $G$ follows from the fact that $x$ and
$y$ are $p^i$-th powers in $K_1$ modulo $K'_1$, and $G$ is of exponent
$p^n$. Thus, if $g\in G$, then
\[[x^{p^{n-i}},g] = [x,g]^{p^{n-i}} = [r^{p^i},g]^{p^{n-i}}=
[r,g]^{p^n} = [r,g^{p^n}] = [r,e] = e\] and analogous with $y$. Since
$K_1$ is generated by $G$, $r$, and $s$, it suffices to show that
$x^{p^{n-i}}$ and $y^{p^{n-i}}$ commute with both $r$ and $s$. That
$x^{p^{n-i}}$ commutes with $r$ is obvious. To verify that it commutes
with $s$ we note that:
\[[x^{p^{n-i}},s] = [(r^{p^i})^{p^{n-i}},s] = [r^{p^n},s] = [r,s]^{p^n} = e\]
since we moded out by $\langle[r,s]^{p^n}\rangle$. An analogous calculation holds
for $y^{p^{n-i}}$.

STEP 4. We adjoin central $p^{n}$-th roots to both $x^{p^{n-i}}$ and
$y^{p^{n-i}}$. We use Lemma \ref{centralroots} and adjoin roots, which
we call $t$ and $v$, respectively. Call the resulting group
$K_2$. In $K_2$ we still have $[r,s]^{p^i}$ in the dominion of $G$ but
not in~$G$, as can be easily verified.

STEP 5. We make $t^{p^i}$ and $v^{p^i}$ into commutators. We can
use Lemma~\ref{makecentralacomm}, since $t$ and $v$ are
central. Denote the resulting group by $K_3$, and write
$t^{p^i}=[q_1,q_2]$, $v^{p^i}=[q_3,q_4]$. Moreover, since
$(t^{p^i})^{p^{n}} = (x^{p^{n-i}})^{p^i} = x^{p^n} = e$, and likewise
$(s^{p^i})^{p^n}=e$, we can choose the $q_i$ to be of exponent $p^n$
each. Again, we still have $[r,s]^{p^i}$ in the dominion of $G$ but
not in $G$. Note, moreover, that $t$ and $v$ are still central in
$K_3$. 

STEP 6. Let $K_4$ be the subgroup of $K_3$ generated by $G$,
$rt^{-1}$, $sv^{-1}$, $q_1$, $q_2$, $q_3$, and $q_4$. 

Note that $rt^{-1}$ and $sv^{-1}$ are of exponent $p^n$. Indeed, since
$t$ and $v$ are central, we have:
\[(rt^{-1})^{p^n} = x^{p^{n-i}}x^{-p^{n-i}}=e=y^{p^{n-i}}y^{-p^{n-i}}=(sv^{-1})^{p^n}.\]
Since each $q_i$ and each element of $G$ is also of
exponent $p^n$, $K_4$ lies in ${\cal N}_2\cap {\cal B}_{p^n}$;
moreover, it is an overgroup of $G$, and
\begin{eqnarray*}
(rt^{-1})^{p^i}[q_1,q_2] &=& r^{p^{i}}t^{-p^i}[q_1,q_2]=r^{p^i}=x,\\
(sv^{-1})^{p^i}[q_3,q_4] & = & s^{p^i}v^{-p^i}[q_3,q_4]= s^{p^i} = y.
\end{eqnarray*}
So $[rt^{-1},sv^{-1}]^{p^i}$ lies in ${\rm dom}_{K_4}^{{\cal
N}_2}(G)$. Since $t$ and $v$ are central in $K_4$,
$[rt^{-1},sv^{-1}]^{p^i} = [r,s]^{p^i}$, which we know does not lie in
$G$; therefore, $G$ cannot be absolutely closed in ${\cal
N}_2\cap{\cal B}_{p^n}$.

This establishes necessity, and proves the theorem.
\end{proof}

We note in particular that in the case $n=1$ the conditions are all
vacuous. Indeed,

\begin{corb}[see Corollary 9.90 in \cite{nildoms}]
If $p$ is an odd prime, then every $G\in {\cal N}_2\cap {\cal B}_p$
is absolutely closed (in that variety).
\end{corb}

\begin{remark}
This corollary also follows from Corollary 1.3 in~\cite{nilexpp}.
\end{remark}

\section{Special bases and absolute closures}

Given the similarity between the characterization of special
amalgamation bases in ${\cal N}_2$ and in ${\cal N}_2\cap{\cal
B}_{p^n}$, a natural question is whether there are any groups which
are special amalgamation bases in the latter but not the former
variety. Note that every special amalgamation base for ${\cal N}_2$
which lies in ${\cal N}_2\cap{\cal B}_{p^n}$ is necessarily also a
special amalgamation base for ${\cal N}_2\cap{\cal B}_{p^n}$. 

When we start looking for possible counterexamples to the converse
implication, we might be tempted to try with easy groups, such as
abelian groups or groups of exponent $p$. However, no such example
will work. Recall that:

\begin{theorem}[Theorems 3.7, 3.13 and 3.17 in \cite{absclosed}]
\label{easygeneral}
Let $G$ be a nilpotent group of class two. 
\begin{itemize}
\item[(i)]\ If $G$ is abelian, then $G$ is absolutely closed in ${\cal
N}_2$ if and only if $G/pG$ is cyclic or trivial for each prime
$p$. Therefore, cyclic groups are absolutely closed in~${\cal N}_2$.
\item[(ii)]\ If $G$ is of exponent $p$, for $p$ a prime, then $G$
is absolutely closed in ${\cal N}_2$ if and only if $Z(G)/G'$ is
cyclic or trivial.
\end{itemize}
\end{theorem}

For $n>1$ we get the same result here:
\begin{theorem}
\label{easyforpn}
Let $G\in{\cal N}_2\cap{\cal B}_{p^n}$, with $p$ an odd prime and
$n>1$. 
\begin{itemize}
\item[(i)]\ If $G$ is abelian, then $G$ is absolutely closed in ${\cal
N}_2\cap{\cal B}_{p^n}$ if and only if $G/pG$ is cyclic or
trivial. Therefore, $G$ is absolutely closed if and only if it is
cyclic or trivial.
\item[(ii)]\ If $G$ is of exponent $p$, then $G$ is absolutely closed in
${\cal N}_2\cap{\cal B}_{p^n}$ if and only if $Z(G)/G'$ is cyclic or
trivial.
\end{itemize}
\end{theorem}

\begin{remark}
Note that if $G$ is of exponent a power of $p$, then $G$ is
$q$-divisible for every integer $q$ relatively prime to $p$;
therefore, if $q$ is a prime different from $p$, then $G/qG$ is
trivial, which means that condition \textit{(i)\/} of
Theorem~\ref{easyforpn} is the same as condition \textit{(i)\/} in
Theorem~\ref{easygeneral}.
\end{remark}

\begin{proof}
\textit{(i)\/} As per the remark above, if $G/pG$ is cyclic or trivial, then $G$
is absolutely closed in ${\cal N}_2$, and hence also in ${\cal
N}_2\cap{\cal B}_{p^n}$. If $G/pG$ is not cyclic, then let $x$ and $y$
be elements of $G$ which project to distinct cyclic summands of
$G/pG$. Then note that a product $x^ay^b$ is a $p$-th power in $G$ if
and only if $p|a$ and $p|b$. Suppose $G$ is absolutely closed. Since
$G$ is abelian, $G$, $x$, $y$, and $p$ fail conditions (a)-(c) of
Theorem \ref{classifabsclosed}. Hence, it must satisfy (d), so there
exist integers $a$, $b$, and $c$, and elements $g_1,g_2\in G$ such that
\begin{eqnarray*}
g_1^p &=& x^a y^b\\
g_2^p &=& x^{b+1} y^c.
\end{eqnarray*}
But that means that $p|b$ and $p|(b+1)$ which is clearly
impossible. So $G$ is not absolutely closed. That this implies that
$G$ is cyclic follows because $G/p^nG$ is a sum of cyclic groups, and
has the same number of summands as $G/pG$, but since $G$ is of
exponent $p^n$ the former group is equal to $G$ itself.

\textit{(ii)\/} Again, if $Z(G)/G'$ is cyclic or trivial, then $G$ is
absolutely closed in ${\cal N}_2$ and hence also in ${\cal
N}_2\cap{\cal B}_{p^n}$. If not, let $x$ and $y$ be elements of $Z(G)$
which project to distinct cyclic generators of $Z(G)/G'$ (which is an
abelian group of exponent $p$). Again, $G$, $x$, $y$ and $p$ do not
satisfy conditions (a)--(c) of Theorem \ref{classifabsclosed}, and the
same argument as above (with congruences rather than equalities)
yields a contradiction to the assumption that it satisfies~(d).
\end{proof}

Fortunately, once we get past the urge to check those easy cases, we find
the examples we were looking for. We need a technical lemma:

\begin{lemmab}[Perturbation argument]
Let $G\in {\cal N}_2$, $H$ a subgroup of $G$, and let $x^q,y^q\in H[G,G]$.
If $h_1, h_1\in H$ then
\[ [x,y]^q\in H \Longleftrightarrow [xh_1,yh_2]^q\in H.\]
\end{lemmab}

\begin{proof}
Note that if $x^q$ and $y^q$ lie in $H[G,G]$, then so do $(xh_1)^q$ and
$(yh_2)^q$, so both commutator brackets lie in the dominion of $H$.

Expanding the bracket bilinearly, we have
\begin{eqnarray*}
[xh_1,yh_2]^q & = & [x,y]^q[h_1,y]^q[x,h_2]^q[h_1,h_2]^q\\
 & = & [x,y]^q[h_1,y^q][x^q,h_2][h_1,h_2]^q.
\end{eqnarray*}
Since $x^q,y^q\in H[G,G]$, and $h_i\in H$, the last three terms on the
right hand side lie in $H$, so the left hand side lies in $H$ if and
only if $[x,y]^q$ lies in $H$, as claimed.
\end{proof}

\begin{theorem}
\label{cantgoup}
Let $p$ be an odd prime. For every $n>0$ there exists a group
$G\in{\cal N}_2\cap{\cal B}_{p^n}$, with $p$ an odd prime, such that
$G$ is absolutely closed in ${\cal N}_2\cap{\cal B}_{p^n}$, but not
absolutely closed in ${\cal N}_2\cap {\cal B}_{p^{n+1}}$. Namely, we
let
\[G = \left\langle x,y \,\bigm|\,
x^{p^n}=y^{p^n}=[x,y]^{p^{n-1}}=e\right\rangle.\]
\end{theorem}

\begin{proof} 
If $n=1$, the result is easy: the group described is a special
amalgamation base of ${\cal N}_2\cap{\cal B}_{p}$ (since everything
there is a special amalgamation base), but not a special amalgamation
base in ${\cal N}_2$ or ${\cal N}_2\cap{\cal B}_{p^m}$ for $m>1$, as it
is abelian but not cyclic. So we may assume that $n>1$.

Consider the group 
\[A=\left\langle a,b \mid a^{p^{n+1}}=b^{p^{n+1}}=[a,b]^{p^{n+1}}=e\right\rangle.\]

Then $G$ is isomorphic to the subgroup generated by $a^p$ and $b^p$;
however, $[a,b]^p\in{\rm dom}_A^{{\cal N}_2}(G)\setminus G$, so $G$
cannot be absolutely closed in ${\cal N}_2\cap{\cal B}_{p^{n+1}}$. We
only need to show now that $G$ is absolutely closed in ${\cal
N}_2\cap{\cal B}_{p^n}$.

Let $K$ be an ${\cal N}_2\cap{\cal B}_{p^n}$-overgroup of $G$, and
suppose that $r^{p^i},s^{p^i}\in G[K,K]$ for some $r,s,\in K$, $1\leq
i \leq n-1$. We may assume that:
\begin{eqnarray*}
r^{p^i} & \equiv & x^{\alpha} y^{\beta} \pmod{K'} \\
s^{p^i} & \equiv & x^{\gamma} y^{\delta} \pmod{K'}.
\end{eqnarray*}
We want to show that $[r,s]^{p^i}$ lies in~$G$. We write
$\mathbf{x}=x^{\alpha} y^{\beta}$ and $\mathbf{y}=x^{\gamma}y^{\delta}$.

Since $r^{p^n}=s^{p^n}=e$, we must have that
$x^{\alpha p^{n-i}}y^{\beta p^{n-i}}$ is central in $G$, and likewise
with $x^{\gamma p^{n-i}}y^{\delta p^{n-i}}$. By calculating the
commutators with $x$ and $y$, we get that $p^{n-1}$ divides $\alpha
p^{n-i}$, $\beta p^{n-i}$, $\gamma p^{n-i}$, and $\delta
p^{n-i}$.
Therefore $p^{i-1}$ divides $\alpha$, $\beta$, $\gamma$, and $\delta$.

So we rewrite the elements
as $\mathbf{x}=x^{\zeta p^{i-1}} y^{\eta p^{i-1}}$,
$\mathbf{y}=x^{\vartheta p^{i-1}}y^{\lambda p^{i-1}}$. Since we can perturb $r$ and
$s$ with elements of $G$ by the Perturbation Argument, we can perturb
$\mathbf{x}$ and~$\mathbf{y}$ by $p^i$-th powers of $x$ and $y$, so we may assume
that $0\leq \zeta,\eta,\vartheta,\lambda < p$.

If $2i\leq n$, then we have that
\[[\mathbf{x},\mathbf{y}]^{p^{n-2i}} = [r^{p^i},s^{p^i}]^{p^{n-2i}} =
[r,s]^{p^n}=e\]
since $K$ is of exponent $p^n$. Substituting the values of
$\mathbf{x}$ and $\mathbf{y}$, we have
\begin{eqnarray*}
[\mathbf{x},\mathbf{y}]^{n-2i} & = & [x^{\zeta p^{i-1}}y^{\eta
p^{i-1}},x^{\vartheta p^{i-1}}y^{\lambda p^{i-i}}]^{p^{n-2i}}\\
& = & [x,y]^{(\zeta\lambda - \eta\vartheta)p^{n-2i+2i-2}}\\
& = & [x,y]^{(\zeta\lambda - \eta\vartheta)p^{n-2}}.
\end{eqnarray*}
Therefore, $p|\zeta\lambda - \eta\vartheta$. If we consider the
vectors $(\zeta,\eta)$ and $(\vartheta,\lambda)$ as being vectors over
$\mathbf{Z}/p\mathbf{Z}$, this means that the two vectors are
proportional. If one of them is the zero vector, then either
$\mathbf{x}$ or $\mathbf{y}$ are the identity element, and there is
nothing to do, e{.}g{.}:
\[ [r,s]^{p^i} = [\mathbf{x},s] = [e,s] = e \in G.\]
So by perturbing $r$ and
$s$ by elements of $G$ we may assume that
$\mathbf{x}^k\equiv\mathbf{y} \pmod{K'}$ for some $k\geq 0$. But
that means that $[r,s]^{p^i}$ actually lies in the dominion of the
subgroup generated by $\mathbf{x}$, which is cyclic and hence
absolutely closed; thus,
$[r,s]^{p^{i}}\in\langle\mathbf{x}\rangle\subseteq G$, as desired.
(In essence, since $\mathbf{y}$ is a power of $\mathbf{x}$, we can
always find exponents to satisfy condition (d); see Theorem 3.7 in
\cite{absclosed}).

If, on the other hand, $2i>n$, then by choosing $j=p^{n-i}$,
$a=b=c=0$, $g_1=e$, and $g_2=x^{\zeta p^{n-i-1}}y^{\eta p^{n-i-1}}$, we
have:
\begin{eqnarray*}
e^{p^i} &\equiv& \mathbf{x}^0\mathbf{y}^0\\
& \equiv & \left(x^{\zeta p^{i-1}}y^{\eta p^{i-1}}\right)^0
\left(x^{\vartheta p^{i-1}}y^{\lambda p^{i-1}}\right)^0 \pmod{G'}\\
\\
\left(x^{\zeta p^{n-i-1}}y^{\eta p^{n-i-1}}\right)^{p^i} & \equiv &
\mathbf{x}^{0+p^{n-i}}\mathbf{y}^0\\
& \equiv &\left(x^{\zeta p^{i-1}}y^{\eta p^{i-1}}\right)^{p^{n-i}} \left(
 x^{\vartheta p^{i-1}}y^{\lambda p^{i-1}}\right)^0 \pmod{G'}
\end{eqnarray*}
so by Lemma~\ref{keylemma},
$[r,s]^{p^{n-i}p^i}=[e,\mathbf{x}][x^{\zeta p^{n-i-1}}y^{\eta
p^{n-i-1}},\mathbf{y}]$. However, \[[r,s]^{p^{n-i}p^i}=[r,s]^{p^n}=e,\]
so we have:
\begin{eqnarray*}
e & = & [r,s]^{p^n}\\
& = & [x^{\zeta p^{n-i-1}}y^{\eta p^{n-i-1}},\mathbf{y}]\\
& = & [x^{\zeta p^{n-i-1}}y^{\eta p^{n-i-1}},x^{\vartheta
p^{i-1}}y^{\lambda p^{i-1}}]\\
& = & [x,y]^{(\zeta\lambda - \eta\vartheta)p^{n-i-1+i-1}}\\
& = & [x,y]^{(\zeta\lambda - \eta\vartheta)p^{n-2}}
\end{eqnarray*}
so again we have that $p|\zeta\lambda - \eta\vartheta$, and we proceed
as we did in the case when $2i\leq n$.

Therefore, $G$ is absolutely closed in ${\cal N}_2\cap{\cal B}_{p^n}$,
as claimed.
\end{proof}

We pause here to adress two questions concerning dominions.

In \cite{episdomssemi}, P.M.~Higgins asks whether the morphic image of
an absolutely closed semigroup is necessarily absolutely
closed. Higgins' question is in the context of the category of all
semigroups, so we cannot answer his question, but we can answer the
corresponding question for arbitrary categories of algebras, and even
for some varieties of semigroups. The answer is negative.

Namely, the variety ${\cal N}_2\cap{\cal B}_{p^n}$ can be defined
entirely with semigroup identities; i{.}e{.} identities in which no
negative number appears as an exponent. To do this, take the group
identities that define the variety, replace $x^{-1}$ with $x^{p^n-1}$
wherever it occurs, and add identities that specify that $x^{p^n}$
acts as a two sided identity. The semigroups which lie in this variety
are just the groups in ${\cal N}_2\cap{\cal B}_{p^n}$, so the
description of dominions in this variety of semigroups is exactly the
same as the description in the variety of groups. In particular, the
description of which groups are absolutely closed is the same as the
one given.

Now consider the ${\cal N}_2$-group presented by
\[G = \langle x,y,z | x^p=y^p=z^p=[x,y]^p=[x,z]=[y,z]=e\rangle.\]
This group is absolutely closed in ${\cal N}_2$ and so also in ${\cal
N}_2\cap{\cal B}_{p^n}$ for any $n$, since it is of exponent $p$ and
$Z(G)/G'$ is generated by the image of $z$. However, its
abelianization is isomorphic to the product of three copies of $Z/pZ$,
and therefore cannot be absolutely closed. Thus, $G$ is absolutely
closed and has a morphic image which is not absolutely closed.

It may be worth noting as well that the subgroup we are moding out by
when considering the morphic image is the commutator subgroup, which
is a fully invariant central subgroup of $G$, and will be normal (in
fact, central) in any overgroup $K$ of $G$ (provided $K\in{\cal
N}_2$). Since dominions respect quotients, any witness to the
fact that $G^{{\rm ab}}$ is not absolutely closed cannot come from an
overgroup $K$ of~$G$. 

Note also that since dominions respect quotients, if $G^{{\rm ab}}$ is
absolutely closed, then $G$ is absolutely closed as well.

We also take this opportunity to address the notion of absolute
closures. Recall that Isbell proves in \cite{isbellone} that in a
right-closed category of algebras (in the sense of Universal Algebra),
every algebra can be embedded in an absolutely closed algebra; since
varieties are right-closed, this implies that any group $G\in{\cal
N}_2\cap{\cal B}_{p^n}$ can be embedded in a group which is absolutely
closed in the same variety (the same holds for ${\cal N}_2$). Isbell
then gives the following definition:

\begin{defb}
An \textit{absolute closure} for the algebra $B$ (in the category
${\cal C}$) is an absolutely closed ${\cal C}$-algebra $D$ containing
$B$, and such that there exists some overalgebra $A$ of $B$, with
$A\in{\cal C}$ and ${\rm dom}_A^{{\cal C}}(B)=D$.
\end{defb}

In other words, an absolute closure for $B$ is an absolutely closed
algebra $D$ which is dominated by $B$ in some overalgebra. In
\cite{isbellone}, Isbell notes that he did not know if absolute
closures existed in general. Soon after, however, examples of algebras
with no absolute closure were constructed.

On the other hand, it is not hard to construct examples of algebras
with an absolute closure, aside from the absolutely closed algebras
themselves. If we work in the variety of all semigroups, then every
group is absolutely closed (this is an easy consequence of the Zigzag
Lemma \cite{isbellone}). If we take a subsemigroup $S$ of a group
$G$, with $S$ not a group, then the dominion of $S$ in $G$ will be the
subgroup generated by $S$; this is absolutely closed (since it is a
group), and properly contains $S$. Thus the subgroup of $G$ generated
by $S$ is an absolute closure for $S$. For example ${\rm
dom}_{\mathbf{Z}}({\mathbf N})=\mathbf{Z}$, and $\mathbf{Z}$ is
absolutely closed, hence an absolute closure for $\mathbf{N}$.

I was unaware that examples of algebras with no absolute closure had
already been constructed, so I began to see if the characterization of
absolutely closed groups in ${\cal N}_2$ given in \cite{absclosed}
could be used to construct such groups and to describe absolute
closures. It turns out, as we will prove below, that in this
variety we have the worst possible case of a negative answer to
Isbell's original question; namely, a group has an absolute closure in
${\cal N}_2$ if and only if it is already absolutely closed. The same
result holds for ${\cal N}_2\cap{\cal B}_{p^n}$.

\begin{theorem}
\label{absclosuresnone}
Let $G\in \mathcal{N}_2$. Then $G$ has an absolute closure in
$\mathcal{N}_2$ if and only if $G$ is already absolutely closed in
$\mathcal{N}_2$. The same is true if we replace $\mathcal{N}_2$ with
${\cal N}_2\cap{\cal B}_{p^k}$ for an odd prime $p$ and $k>0$.
\end{theorem}

\begin{proof} We prove the ${\cal N}_2$ case first.
One implication is trivial. So suppose that $G$ is not absolutely
closed, and let $K$ be any nil-2 overgroup of~$G$. Let $D_K={\rm
dom}_K^{{\mathcal{N}}_2}(G)$. We want to show that $D_K$ is also not
absolutely closed.

If $G=D_K$, then there is nothing to prove. So we may assume that $G$
is properly contained in $D_K$. By the description of dominions from
Theorem \ref{nildomsdesc}, there exist elements $r,s\in K$, and $n>0$
such that $r^n$ and $s^n$ lie in $G[K,K]$, and $[r,s]^n$ does not lie
in~$G$ (it lies in $D_K$, however). We write $[r,s]^n=d_0$. Note that
the description of dominions also implies that $[D_K,D_K]=[G,G]$.

Choose $r',s'\in [K,K]$ such that $r^nr',s^ns'\in G$, and write
$r^nr'=x$, $s^ns'=y$. We claim that $x$, $y$, and $n$ do not satisfy
either condition~\textit{(i)\/} or condition~\textit{(ii)\/} from Theorem
\ref{classiffullvar} with respect to $D_K$. This will show that $D_K$
is not absolutely closed, as these conditions describe the nil-2
absolutely closed groups.

Condition \textit{(i)\/} is the easiest to handle: condition
\textit{(i)\/} is equivalent to the statement that no nil-2 extension
of~$D_K$ contains $n$-th roots for both $x$ and~$y$. However, $x$ and
$y$ are both $n$-th powers modulo the commutator in~$K$, so there is a
nil-2 extension of~$K$ with $n$-th roots for both $x$ and $y$.  In
particular, there is such an extension of $D_K$ as well. So $x$, $y$,
and $n$ do not satisfy condition \textit{(i).}

Condition \textit{(ii)\/} is somewhat more difficult, and we proceed by
contradiction. Suppose that condition \textit{(ii)\/} is satisfied
in~$D_K$. Then there exist integers $a$, $b$, and~$c$, and elements
$d_1, d_2\in D_K$ such that
\begin{eqnarray*}
{d_1}^n & \equiv & x^a y^b \pmod{{D_K}'}\\
{d_2}^n & \equiv & x^{b+1}y^c \pmod{{D_K}'}.\\
\end{eqnarray*}
By Theorem \ref{nildomsdesc}, we can write $d_i$ as
\[
d_i = g_i[k_{i11}, k_{i12}]^{n_{i1}}\cdots
[k_{is_i1},k_{is_i2}]^{n_{is_i}}\]
for some $g_i\in G$, $s_i\geq 0$, $k_{ijk}\in K$, with
$k_{ijk}^{n_{ij}}\in G[K,K]$. In particular, $d_i\equiv g_i\pmod{K'}$.

By Lemma \ref{keylemma}, the congruences above imply that
$[r,s]^n=[d_1,x][d_2,y]$. However, since $d_i\equiv g_i\pmod{K'}$,
this means that $[r,s]^n=[g_1,x][g_2,y]$. But the right hand side
clearly lies in~$G$, and this contradicts the choice of $r$ and
$s$. So $D_K$, $x$, $y$, and $n$ cannot satisfy (ii) either.

Therefore $D_K$ is not absolutely closed. We conclude that
$G$ has no absolute closure, which proves the ${\cal N}_2$ case.

For the ${\cal N}_2\cap{\cal B}_{p^k}$ case, we proceed exactly as
above, replacing $n$ with $p^i$, and restricting $K$ to the variety in
question. The argument above shows that $G$, $x$, $y$ and $p^i$ do not
satisfy (a) nor (d) of Theorem \ref{classifabsclosed}. If
$[x,y]^{p^{\alpha}}\not=e$, then $[r,s]^{p^{2i+\alpha}}\not=e$, which
contradicts the fact that $K$ is of exponent $p^k$, since $k\leq
2i+\alpha$. Thus, they cannot satisfy (b). For (c), if $2i>k$ and
there is a solution to the congruences with $[d_1,x][d_2,y]\not=e$,
then by Lemma~\ref{keylemma} we have that $[r,s]^{jp^i}\not=e$. But
$p^k|jp^i$, which again contradicts the fact that $K$ is of exponent
$p^k$. Therefore, $D_K$ is also not absolutely closed in this case,
and this proves the theorem in full.
\end{proof}

\section{Weak and strong bases}

The case of weak and strong amalgamation bases turns out to be
remarkably simpler than the special amalgamation case. The results are
just simple applications of our observations at the beginning, plus
some trivial checks on the known conditions for embeddability. We
include them for completeness, and to note some differences between
this case and the special amalgamation case above.

Necessary and sufficient conditions for strong embeddability of an
amalgam of ${\cal N}_2$-groups were first given by Wiegold in
\cite{wiegoldamalg}. However, they were hard to work with, as they
involved checking the existence of maps with certain properties on
some tensor products. Later, B.~Maier gave easier conditions for
both weak and strong embeddability. We recall those conditions now:

\begin{theorem}[B. Maier, Hauptsatz in \cite{amalgone}]
\label{maierweak}
The amalgam $(G,K;H)$ of nil-$2$ groups is weakly embeddable in ${\cal
N}_2$ if and only if
\begin{itemize}
\item[(i)]\ $G'\cap H\subseteq Z(K)$ and $K'\cap H\subseteq Z(G)$; and
\item[(ii)]\ For all $s>0$, $q_i>0$, $g_i\in G$, $g_i'\in G'$, $k_i\in
K$, $k_i'\in K'$ with $g_i^{q_i}g_i', k_i^{q_i}k_i'\in H$, $1\leq
i\leq s$, for all $h\in H$ we have
\[\prod_{i=1}^s [g_i,k_i^{q_i}k_i'] = h \Longleftrightarrow \prod_{i=1}^s
[g_i^{q_i}g_i',k_i]=h.\]
\end{itemize}
\end{theorem}

\begin{theorem}[B. Maier, Satz 3 in \cite{amalgtwo}]
\label{maierstrong}
The amalgam $(G,K;H)$ of nil-2 groups is strongly embeddable in ${\cal N}_2$ if and
only if
\begin{itemize}
\item[(i)]\ $G'\cap H\subseteq Z(K)$ and $K'\cap H\subseteq Z(G)$;
\item[(ii)]\ For each prime power $p^i$, if $g\in G$, $g'\in G'$, $k\in
K$, $k'\in K'$ are such that $g^{p^i}g', k^{p^i}k'\in H$, we have
\[[g^{p^i}g',k] = [g,k^{p^i}k'] \in H.\]
\end{itemize}
\end{theorem}

We can make some simplifying observations:

\begin{remark}
\label{simplifweak} 
In Theorem~\ref{maierweak}, if $H$ is $k$-divisible, then we can
restrict $q_i$ to be different from $k$. This because if $H$ is $q_i$
divisible, we can replace $g^{q_i}g'$ and $k^{q_i}k'$ with some
$q_i$-th power of an element of $H$, and then that commutator will
certainly lie in $H$, and we can omit it from the given condition. In
particular, if $H$ is (or if $G$ and $K$ are) of exponent $n$, we can
restrict $q_i$ to numbers not relatively prime to $n$. If $G$ and $K$
both lie in ${\cal B}_{p^n}$, we can restrict the $q_i$ to be powers
of $p$. In that case, we can further restrict them to being powers
$p^i$ with $i<n$, since otherwise the corresponding commutator is
trivial.
\end{remark}

\begin{remark}
\label{simplifstrong} 
Likewise, in Theorem \ref{maierstrong}, if $G$ and $K$ are both of
finite exponent $n$, we can restrict the primes $p$ to prime factors
of $n$. If $G$ and $K$ are of exponent $p^n$, we can also restrict the
prime powers to $p^i$, with $1\leq i\leq n-1$.
\end{remark}

\begin{remark}
The conditions and simplifications given reduce, in the case $n=1$, to
the existence of a ``coupled central series'' of $G$ and $K$ over $H$,
as described in Definition 1.1 and Theorem 1.2 in \cite{nilexpp}.
\end{remark}

The characterization of weak and strong bases in ${\cal N}_2$ was
given by Saracino:

\begin{theorem}[Saracino, Theorem 3.3 in \cite{saracino}]
\label{saracinobases}
Let $G\in {\cal N}_2$. The following are equivalent:
\begin{itemize}
\item[(a)]\ $G$ is a weak amalgamation base for ${\cal N}_2$.
\item[(b)]\ $G$ satisfies $G'=Z(G)$, and $\forall g\in G, \forall n>0$
($g\in G^nG'$ or $\exists y\in G$ and $\exists k\in\mathbf{Z}$ such
that ($y^n\equiv g^k \pmod{G'}$ and $[y,g]\not=e$)).
\item[(c)]\ $G$ satisfies $G'=Z(G)$, and for all $g\in G$ and $n>0$,
either $g$ has an $n$-th root modulo $G'$ or else $g$ has no $n$-th
root in any overgroup $K\in{\cal N}_2$ of $G$.
\item[(d)]\ $G$ is a strong amalgamation base for ${\cal N}_2$.
\end{itemize}
\end{theorem}

\begin{remark}
Again it is not hard to verify that in Theorem \ref{saracinobases} we
can restrict $n$ to prime powers. Moreover, the conditions for
adjunction of roots show that if $G$ is of exponent $n$, then an
element of $G$ has an $n$-th root in some ${\cal N}_2$-overgroup if
and only if it is central. An easy way to verify it is to note that an
element $x$ which is a $k$-th power in some extension must commute with
every element of exponent $k$ in $G$; for if $r$ is such a root in
some overgroup, and $y$ is of exponent $k$, then in that overgroup we have
\[[x,y] = [r^k,y] = [r,y^k] = [r,e] = e.\]
In fact, an $x$ which is a $k$-th power in some ${\cal N}_2$-extension
of $G$ must commute with any element whose image in the
abelianization of $G$ has exponent $k$. So if $G$ is of exponent $n$,
any element to which we can adjoin an $n$-th root must be central, and
of course we can always adjoin $n$-th roots to central elements. But
if $G'=Z(G)$, then in such a group an element has an $n$-th root in
some extension if and only if it is a commutator, and therefore will
always lie in $G^nG'=G'$. Thus, for $G$ of exponent $p^k$, we can
restrict the prime powers in Theorem \ref{saracinobases} to $p^i$,
with $0< i < k$. 
\end{remark}

Here we want a description of the amalgamation bases in ${\cal
N}_2\cap{\cal B}_{p^n}$, and we might also ask if there are any which
are strong or weak bases there but not in ${\cal N}_2$.  However, in
contrast to Theorem \ref{cantgoup}, we have:

\begin{theorem}[cf. Theorem \ref{cantgoup}]
\label{weakandstrongaresame}
Let $p$ be an odd prime, $n>0$, and let $G\in{\cal N}_2\cap{\cal
B}_{p^n}$. The following are equivalent:
\begin{itemize}
\item[(i)]\ $G$ is a weak amalgamation base for ${\cal N}_2$.
\item[(ii)]\ $G$ is a weak amalgamation base for ${\cal N}_2\cap{\cal
B}_{p^n}$. 
\item[(iii)]\ $G$ satisfies $Z(G)=G'$ and for each $i$, $1\leq i\leq
n-1$, and each $g\in G$, either $g\in G^{p^i}G'$ or else there exists
$c>0$ and $y\in G$ with $y^{p^i}\equiv g^c\pmod{G'}$ and
$[y,g]\not=e$. 
\item[(iv)]\ $G$ is a strong amalgamation base for ${\cal N}_2\cap
{\cal B}_{p^n}$. 
\item[(v)]\ $G$ is a strong amalgamation base for ${\cal N}_2$.
\end{itemize}
\end{theorem}

\begin{proof}
Clearly $(v)\Rightarrow(i)\Rightarrow(ii)$ and
$(v)\Rightarrow (iv)\Rightarrow(ii)$. We show that
$(ii)\Rightarrow(iii)\Rightarrow(v)$. 

Suppose that $G$ does not satisfy \textit{(iii).} If $Z(G)\not=G'$, then we
take an element $x\in Z(G)\setminus G'$. We let $K_1$ be an
${\cal N}_2\cap{\cal B}_{p^n}$-overgroup in which $x$ is a commutator,
using Lemma \ref{makecentralacomm}, and we let $K_2$ be a ${\cal
N}_2\cap{\cal B}_{p^n}$-overgroup in which $x$ is not central, say
$G\coprod^{{\cal N}_2}Z/p^nZ$. Then the amalgam $(K_1,K_2;G)$ cannot
be embedded in an ${\cal N}_2$-group, since $x$ would have to be central by
virtue of being a commutator in $K_1$, but would have to be noncentral
since it is not central in $K_2$. So $G$ cannot satisfy \textit{(ii).}

Suppose instead that $Z(G)=G'$, and that $g\in G$ is such that
$g\notin G^{p^i}G'$, but there is some ${\cal N}_2$-extension with a
$p^i$-th root for $g$. Note that since there is some extension with a
$p^i$-th root for $g$, $g$ must commute with every element of exponent
$p^i$-th in $G$. Since $G$ is of exponent $p^n$, this means that
$x^{p^{n-i}}$ must be central in $G$, hence lies in the commutator of~$G$.

Let $K_1$ be an overgroup in which $g$ does not commute with an
element of exponent $p^i$, for example, consider $G\coprod^{{\cal
N}_2}Z/p^iZ$. Since $g\notin G^{p^i}G'$, $g$ does not commute with the
generator of the cyclic group. For $K_2$ we want a group of exponent
$p^n$ where $x$ has a $p^i$-th root modulo $K_2'$, since then we will
have a problem: in any group into which we embedd the amalgam
$(K_1,K_2;G)$, $x$ would be a $p^i$-th power modulo the commutator, so
it would have to commute with everything of exponent $p^i$, yet it
does not do so in~$K_1$.

To construct $K_2$ we again have to be a bit careful.  First consider
a nil-2 group in which $x$ has a $p^i$-th root, and call such a
root $r$. We may assume that the group is generated by $G$ and
$r$. Since $x^{p^{n-i}}$ lies in $G'$, it is central, and we can
adjoin a central $p^n$-th root $t$ to that element. Since $G$ is of
exponent $p^n$, $t^{p^i}$ is of exponent $p^n$; using Lemma
\ref{makecentralacomm}, we can make $t^{p^i}$ into a commutator,
$[q_1,q_2]$, with $q_i$ of exponent $p^n$. Now we consider the group
generated by $G$, $q_1$, $q_2$, and $rt^{-1}$. Then
\[(rt^{-1})^{p^n}=  r^{p^n} t^{-p^n}= x^{p^{n-i}}x^{-p^{n-i}} = e.\]
Therefore, this group is of exponent $p^n$, since it is generated by
elements of exponent $p^n$; moreover,
$x=(rt^{-1})^{p^i}[q_1,q_2]$, so $x$ is indeed a $p^i$-th power modulo
the commutator subgroup. Make this group $K_2$, and we are done.

Thus, if $G$ does not satisfy \textit{(iii)}, it cannot be a weak amalgamation
base.

Now suppose that $G$ satisfies \textit{(iii).} We use Theorem \ref{maierstrong}
to show that $G$ is a strong amalgamation base in~${\cal N}_2$. Let $K_1$ and $K_2$ be
overgroups of $G$. Then $K_1'\cap G\subseteq Z(G)=G'$,
hence lies in the center of $K_2$, and likewise for $K_2'\cap
G$. 
Suppose that $r\in K_1$, $r'\in K_1'$, $s\in K_2$, $s'\in K_2'$
are such that $r^{q}r'=x\in G$ and $s^{q}s'=y\in G$ for some prime
power $q$. If $\gcd(p,q)=1$, then since $G$ is $q$-divisible, we can
write $x=g_1^q$ and $y=g_2^q$. Then:
\[[r^qr',s]=[g_1^q,s]=[g_1,s^qs']=[g_1,g_2^q]=[g_1,g_2]^q\in G\]
and likewise with $[r,s^qs']$. So we may assume that $q=p^k$. If
$k\geq n$, then $x$ and $y$ must be central in $G$, and therefore they
lie in $G'$. So $[r^qr',s]=[r,s^qs']=e$. We may therefore assume that
$q=p^i$ with $1\leq i \leq n-1$.

Clearly there are overgroups of $K_1$ in which $x$ has a $p^i$-th
root, so we must have $x\in G^{p^i}G'$. Let $g\in G$ such that
$g^{p^i}\equiv x\pmod{G'}$. Likewise, there is an $h\in G$ such that
$h^{p^i}\equiv y\pmod{G'}$. Then
\[[r^{p^i}r',s] = [x,s] = [g^{p^i},s] = [g,s^{p^i}] = [g,s^{p^i}s'] =
[g,h^p]=[g,h]^p\in G\]
and
\[[r,s^{p^i}s'] = [r,y] = [r,h^{p^i}] = [r^{p^i},h] = [r^{p^i}r',h] =
[g^p,h]= [g,h]^p\in G.\]
Thus the amalgam $(K_1,K_2;G)$ is strongly embeddable into an ${\cal
N}_2$-group. This proves that $G$ is a strong amalgamation base for
${\cal N}_2$, and proves the theorem. 
\end{proof}

\begin{remark}
It is worth noting that for the case $n=1$, the conditions given
reduce to $G'=Z(G)$. This is consistent with Maier's Theorem 2.1 in
\cite{nilexpp} and with Saracino's Proposition 3.4 in \cite{saracino}.
\end{remark}

\begin{corb}
\label{enoughcheckn}
A group $G\in{\cal N}_2\cap{\cal B}_{p^n}$ is a strong amalgamation
base in ${\cal N}_2$ if and only if every amalgam of ${\cal
N}_2\cap{\cal B}_{p^n}$ groups with core $G$ can be (strongly) embedded
into an ${\cal N}_2$-group.
\end{corb}

\section{Comments and questions}

Using Proposition \ref{pparts}, we see that the results in the previous two
sections characterize amalgamation bases (weak, strong, and special)
in the variety ${\cal N}_2\cap{\cal B}_n$ for any odd number
$n$. Since a group of exponent $2$ is necessarily abelian, and in
abelian groups of exponent two every amalgam can be strongly embedded,
the results trivially extend to the case where $n$ is not divisible by
$4$.

In order to extend it to all $n$, it would suffice to address the case
of $n=2^{i+1}$ with $i>0$. Unfortunately, that case turns out to be much
more complicated than the odd exponent case.

To make the arguments easier to state, we will write ``no ${\cal
N}_2$-extension of $G$ has a $2^i$-th root for $x$'' rather than
writing the equivalent condition given by Theorem~\ref{genrootadjunction}.

One problem is clear when we note that in a group of exponent $2^{n+1}$,
the commutators are necessarily of exponent $2^{n}$, and that it is no
longer true that a group generated by elements of order $2^{n+1}$ is
necessarily of order $2^{n+1}$. If we attempt to proceed naively as we did
above, we will get that the conditions under which an amalgam of
exponent $2^{n+1}$ groups can be embedded \textit{into a group of
exponent $2^{n+2}$} are the same as for ${\cal N}_2$, but it is
unclear how to lower that last exponent by half.

Thus, for example, I do not know if the description of dominions is
still true for groups of exponent 4 (in which case, dominions would be
trivial, since squares are central in such a group). 

There are other difficulties. For example,
Corollary~\ref{enoughcheckn} no longer holds for groups of exponent
$2^{n+1}$ (we will show an example below). One difficulty lies in the
condition $Z(G)=G'$. In the case of odd exponent, if $Z(G)\not=G'$, we
took a central element which was not a commutator and extended $G$ in
one direction to make that element a commutator, and in the other
direction to make it non-central. However, if $G$ is of exponent
$2^{n+1}$, there could be a central element of order $2^{n+1}$; but
all commutators are of exponent $2^{n}$, so we would not be able to
extend $G$ to a group of exponent $2^{n+1}$ in which that element is a
commutator.  If we allow the $K_i$ to be of exponent $2^{n+2}$, so that
commutators are of exponent $2^{n+1}$, the difficulty disappears,
although we can now only guarantee that the resulting amalgam is
embeddable in a group of exponent $2^{n+3}$. We will also be able to
use this to handle a failure of the second part of the
condition. Thus, by proceeding as we did above, we obtain the
following result:

\begin{theorem}
Let $G$ be an ${\cal N}_2$-group of exponent $2^{n+1}$, with $n>0$. The
following are equivalent:
\begin{itemize}
\item[(i)]\ $G$ is a weak amalgamation base for ${\cal N}_2$.
\item[(ii)]\ Every amalgam of ${\cal N}_2\cap{\cal B}_{2^{n+2}}$
groups with core $G$ is weakly embeddable in ${\cal N}_2$.
\item[(iii)]\ $G'=Z(G)$ and for every $g\in G$ and every $i$, $1\leq
i\leq n$, either $g\in G^{2^i}G'$ or else no ${\cal N}_2$-extension
of $G$ has a $2^i$-th root for $g$.
\item[(iv)]\ Every amalgam of ${\cal N}_2\cap{\cal B}_{2^{n+2}}$
groups with core $G$ is strongly embeddable in ${\cal N}_2$.
\item[(v)]\ $G$ is a strong amalgamation base for ${\cal N}_2$.
\end{itemize}
\end{theorem} 

Another difficulty lies in trying to characterize the elements $x\in
G$ which can be a $2^i$-th power modulo the commutator in an extension
of $G$ of exponent $2^{n+1}$. Of course, it must be true that there is
some ${\cal N}_2$-extension where $x$ is a $2^i$-th power; but we also
get two extra conditions that are not implied by that statement:
$x^{2^{n}}$ must be trivial (if $i>0$), and $x^{2^{n-i}}$ must be central
in $G$ (in fact, in the extension). 

Although Corollary~\ref{enoughcheckn} no longer holds for power of two
exponent, we can still characterize the groups $G\in{\cal N}_2\cap{\cal
B}_{2^{n+1}}$ for which every amalgam of ${\cal
N}_2\cap{\cal B}_{2^{n+1}}$-groups with core $G$ are embeddable in
${\cal N}_2$. We define
\[\Omega^i(G)=\{g\in G\mid g^{2^i}=e\}.\] 

\begin{theorem}
Let $G\in{\cal N}_2\cap{\cal B}_{2^{n+1}}$, with $n>0$. Then the
following are equivalent:
\begin{itemize}
\item[(i)]\ Every amalgam $(K_1,K_2;G)$ of ${\cal N}_2\cap{\cal
B}_{2^{n+1}}$ groups is weakly embeddable in ${\cal N}_2$.
\item[(ii)]\ Every amalgam $(K_1,K_2;G)$ of ${\cal N}_2\cap{\cal
B}_{2^{n+1}}$ groups is strongly embeddable in ${\cal N}_2$.
\item[(iii)]\ $\Omega^n(Z(G))=G^{2^n}G'$ and for each $x\in G$ and each
$i$, $1\leq i\leq n-1$, one of the following holds:
\begin{itemize}
\item[(a)]\ $x\notin G^{2^i}G'$; or
\item[(b)]\ $x^{2^{n-i}}\notin Z(G)$; or 
\item[(c)]\ $x^{2^n}\not=e$; or 
\item[(d)]\ No ${\cal N}_2$ extension of $G$ has a
$2^i$-th root for $x$.
\end{itemize}
\end{itemize}
\end{theorem}

\begin{proof}
Note that for any $G\in{\cal N}_2\cap{\cal B}_{2^{n+1}}$,
$G^{2^n}G'\subseteq \Omega^n(Z(G))$, as every $2^n$-th power is
central and of exponent $2$, and every commutator is central and of
exponent $2^n$.

Clearly $(ii)\Rightarrow(i)$. Suppose $G$ fails to satisfy
$(iii)$. If $\Omega^n(Z(G))\not=G^{2^n}G'$, let $g\in G$ be central of
exponent $2^n$, but not in $G^{2^n}G'$. Since $g$ is of exponent $2^n$
and central, we can extend $G$ to a group where $g$ is a commutator;
moreover, we can make $g=[q_1,q_2]$ with $q_i$ of exponent $2^n$, so
the resulting group is still of exponent $2^{n+1}$. Let that group be
$K_1$. Then let $K_2=G\coprod^{{\cal N}_2}Z/2^nZ$; $K_2$ is
also of exponent $2^{n+1}$, and $g$ is not central in $K_2$. Thus
$(K_1,K_2;G)$ cannot be weakly embedded.

If $\Omega^n(Z(G))=G^{2^n}G'$, let $x\in G$ be such that $x\notin
G^{2^i}G'$, but such that $x^{2^n}=e$, $x^{2^{n-i}}$ is central in
$G$, and there is some extension of $G$ where $x$ has
a $2^i$-th root.
Let $K_1$ be obtained from $G$ by first adjoining a $2^i$-th root
$r$ to $x$; then adjoining central $2^n$-th root $t$ to $x^{2^{n-i}}$.
Note that $t^{2^i}$ is central of exponent $2^n$, since $t^{2^{i+n}} =
{x^{2^n}}=e$, so we can make $t^{2^i}$ into a commutator $[q_1,q_2]$,
with $q_i$ of exponent ${2^n}$. Then the group generated by $G$,
$rt^{-1}$, $q_1$ and $q_2$ is of exponent $2^{n+1}$ (since $rt^{-1}$
is also of exponent $2^n$), and here $x=r^{2^i}[q_1,q_2]$. Let
$K_2=G\coprod^{{\cal N}_2}Z/2^iZ$; then, in $K_2$, $x$ does not commute
with an element of exponent $2^i$, which gives that $(K_1,K_2;G)$ is
not embeddable in ${\cal N}_2$. This proves that
$(i)\Rightarrow(iii)$.

To prove that $(iii)\Rightarrow (ii)$, let $(K_1,K_2;G)$ be an amalgam
of ${\cal N}_2\cap{\cal B}_{2^{n+1}}$ groups, and assume $G$ satisfies
$(iii)$. We verify Maier's conditions for strong embeddability.

The elements of $K_1'\cap G$ are central, and of exponent $2^n$, so
they lie in $\Omega^n(Z(G))=G^{2^n}G'$. However, $G^{2^n}G'$ is
central in $K_2$, so $K_1'\cap G\subseteq Z(K_2)$. A symmetric
argument yields $K_2'\cap G\subseteq Z(K_1)$.

Now suppose that $k_1\in K_1$, $k_1'\in K_1'$, $k_2\in K_2$, $k_2'\in
K_2'$ and $q>0$ is a prime power such that $k_1^qk_1',k_2^qk_2'\in
G$. Clearly we may assume that $q$ is a power of~2, and moreover
that $q=2^i$ with $1\leq i\leq n-1$ ($2^n$-th powers are central in
$G$ and of exponent at most 2, so they lie in $\Omega^n(Z(G))$, and
so are central in each $K_i$). Let
$x=k_1^qk_1'$, $y=k_2^qk_2'$. It is easy to verify that
$x^{2^n}=y^{2^n}=e$ and that $x^{2^{n-i}}$ and $y^{2^{n-i}}$ are
central in $G$, so we must have $x,y\in G^{2^i}G'$. Now we can proceed
as before to get that $[k_1^qk_1',k_2]=[k_1,k_2^qk_2']\in G$, and so
the amalgam is strongly embeddable in ${\cal N}_2$.
\end{proof}

\begin{remark}
For $n=1$, condition \textit{(ii)\/} reduces to the condition
$\Omega^1(Z(G))=G^2G'$.
\end{remark}

\begin{example}
For example, the group
\[G = \left\langle a,b,c \mid a^4=b^4=c^4=[a,b]^2=[a,c]=[b,c]=e;
c^2=[a,b]\right\rangle \] is such that every amalgam of groups of
exponent $4$ with core $G$ is embeddable (into a group of exponent at
most~8). This because $\Omega^1(Z(G))=\langle
a^2,b^2,[a,b]\rangle$. However, it is clearly not a strong
amalgamation base for ${\cal N}_2$, since $c$ is central but not a
commutator. This example shows that Corollary~\ref{enoughcheckn} does
not hold for power-of-two exponent.
\end{example}

The results above show that we will encounter new difficulties
in trying to extend the results in Sections 2 and 4 to ${\cal
N}_2\cap{\cal B}_{2^{n+1}}$.
We therefore have the following questions:

\begin{openquestion} What are the conditions for weak and strong embeddability
of amalgams of ${\cal N}_2\cap{\cal B}_{2^{n+1}}$ groups?
What is a description of dominions in the variety
${\cal N}_2\cap{\cal B}_{2^{n+1}}$?
What is a characterization the weak, strong, and special
amalgamation bases for the varieties ${\cal N}_2\cap{\cal B}_{2^{n+1}}$?
\end{openquestion}

The varieties we have considered so far are not the only subvarieties
of ${\cal N}_2$. The subvariety of ${\cal N}_2$ are in one-to-one
correspondence with pairs of nonnegative integers $(m,n)$, where
$n|m/{\rm gcd}(2,m)$. To each such pair corresponds the variety
defined by the identities:
\[x^m = [x_1,x_2]^n = [x_1,x_2,x_3] = e\]
(see for example \cite{jonsson} and \cite{classifthree}). So far we
have only considered the cases of $m=n$ with $m$ odd, and $m=n=0$. So we
also ask
\begin{openquestion} 
For the subvarieties of ${\cal N}_2$ we have not considered, find
conditions for weak and strong embeddabilty of amalgams. Describe
dominions in those subvarieties, and characterize the weak, strong,
and special amalgamation bases.
\end{openquestion}

We note that if $n$ is not square free, we already know that there are
nontrivial dominions in such a variety (Theorem 9.92 in~\cite{nildoms}),
so describing dominions and characterizing special amalgamation bases
is a nontrivial problem there.




\begin{references}

\bibitem{episdomssemi}
Peter~M. Higgins, {Epimorphisms, dominions and semigroups}, \emph{Algebra
  Universalis} \textbf{21} (1985), no.~2--3, 225--233, {MR {\bf 87j}:20105}.

\bibitem{episandamalgs}
Peter~M. Higgins, {Epimorphisms and amalgams}, \emph{Colloq. Math.} \textbf{56} (1988),
  no.~1, 1--17, {MR {\bf 89m}:20083}.

\bibitem{isbellone}
J.R. Isbell, {Epimorphisms and dominions}, \emph{in} Proc. of the Conference on
  Categorical Algebra, La Jolla 1965 (S.~Eilenberg et~al, ed.), Lange and
  Springer, New York, 1966, pp.~232--246, {MR {\bf 35}:105a} ({T}he statement of the {Z}igzag
  {L}emma for {\it rings} in this paper is incorrect. {T}he correct version is
  in \cite{isbellfour}).

\bibitem{isbellfour}
J.R. Isbell, {Epimorphisms and dominions {IV}}, \emph{J. London Math. Soc.} (2)
  \textbf{1} (1969), 265--273, {MR {\bf 41}:1774}.

\bibitem{jonsson}
Bjarni J{\'o}nsson, {Varieties of groups of nilpotency three}, \emph{Notices
  Amer. Math. Soc.} \textbf{13} (1966), 488.

\bibitem{machenry}
T.~Mac{H}enry, {The tensor product and the 2nd nilpotent product of
  groups}, \emph{Math. Z.} \textbf{73} (1960), 134--145, {MR {\bf 22}:11027a}.

\bibitem{absclosed}
Arturo Magidin, {Absolutely closed nil-2 groups}, \emph{Algebra Universalis}
  \textbf{42} (1999), no.~1--2, 61--77.

\bibitem{nildoms}
Arturo Magidin, {Dominions in varieties of nilpotent groups}, \emph{Comm. Algebra}
  \textbf{28} (2000), no.~3, 1241--1270.

\bibitem{amalgone}
Berthold~J. Maier, {Amalgame nilpotenter {G}ruppen der {K}lasse zwei},
  \emph{Publ. Math. Debrecen} \textbf{31} (1985), 57--70, {MR {\bf 85k}:20117}.

\bibitem{amalgtwo}
Berthold~J. Maier, {Amalgame nilpotenter {G}ruppen der {K}lasse zwei {II}}, \emph{Publ.
  Math. Debrecen} \textbf{33} (1986), 43--52, {MR {\bf 87k}:20050}.

\bibitem{nilexpp}
Berthold~J. Maier, {On nilpotent groups of exponent {$p$}}, \emph{J. Algebra} \textbf{127}
  (1989), 279--289, {MR {\bf 91b}:20046}.

\bibitem{nonembed}
B.H. Neumann and James Wiegold, {On certain embeddability criteria for
  group amalgams}, \emph{Pub. Math. Debrecen} \textbf{9} (1962), 57--64, {MR {\bf
  26}:3763}.

\bibitem{classifthree}
V.N. Remeslennikov, {Two remarks on 3-step nilpotent groups}, \emph{Algebra i
  Logika} (1965), no.~2, 59--65 (Russian), {MR {\bf 31}:4838}.

\bibitem{saracino}
D.~Saracino, {Amalgamation bases for nil{-$2$} groups}, \emph{Algebra
  Universalis} \textbf{16} (1983), 47--62, {MR {\bf 84i}:20035}.

\bibitem{wiegoldamalg}
James Wiegold, {Nilpotent products of groups with amalgamations}, \emph{Publ.
  Math. Debrecen} \textbf{6} (1959), 131--168, {MR {\bf 21}:3478}.

\end{references}

\nocite{isbellfour}
\bibliographystyle{amsplain}

\end{article}
\end{document}